\documentclass[12pt]{article}

\usepackage{amsmath,amsthm}
\usepackage{amssymb}
\usepackage{bm}
\usepackage{bbm}
\usepackage{graphicx}
\AtBeginDocument{
        \DeclareGraphicsExtensions{.pdf,.jpg,.png}}

\usepackage{enumitem}

\usepackage{hyperref}

\usepackage[T1]{fontenc}

\usepackage{listings}
\usepackage{color} 
\definecolor{mygreen}{RGB}{28,172,0} 
\definecolor{mylilas}{RGB}{170,55,241}

\newtheorem{theorem}{Theorem}[section]

\newtheorem{lemma}[theorem]{Lemma}
\newtheorem{prop}[theorem]{Proposition}

\theoremstyle{definition}
\newtheorem{definition}[theorem]{Definition}

\numberwithin{equation}{section}

\newcommand{\R}{\mathbb{R}}

\newcommand{\C}{\mathbb{C}}

\newcommand {\aplt} {\ {\raise-.5ex\hbox{$\buildrel<\over\sim$}}\ } 
\def\dddots{\mathinner{\mkern1mu\raise\p@  
    \hbox{.}\mkern2mu\raise4\p@\hbox{.}\mkern2mu
    \raise7\p@\vbox{\kern7\p@\hbox{.}}\mkern1mu}}%

\begin{document}


\baselineskip=17pt


\title{On the Growth of Lebesgue Constants for Degree One Fekete Points in terms of the Dimension}

\author{Len Bos\\
Department of Computer Science\\ 
Univesity of Verona\\
Italy\\
\\
E-mail: leonardpeter.bos@univr.it}
\maketitle
\date{}

\renewcommand{\thefootnote}{}

\footnote{2010 \emph{Mathematics Subject Classification}: Primary 41A17; Secondary 41A63.}

\footnote{\emph{Key words and phrases}: Fekete points, optimal measures,  optimal experimental design,  simplex.}

\renewcommand{\thefootnote}{\arabic{footnote}}
\setcounter{footnote}{0}

\lstset{language=Matlab,%
    breaklines=true,%
    morekeywords={matlab2tikz},
    keywordstyle=\color{blue},%
    morekeywords=[2]{1}, keywordstyle=[2]{\color{black}},
    identifierstyle=\color{black},%
    stringstyle=\color{mylilas},
    commentstyle=\color{mygreen},%
    showstringspaces=false,
    numbers=left,%
    numberstyle={\tiny \color{black}},
    numbersep=9pt, 
    emph=[1]{for,end,break},emphstyle=[1]\color{red}, 
}


\begin{abstract}
We discuss the growth of the  Lebesgue constants for polynomial interpolation at Fekete points for fixed degree (one) and varying dimension, and underlying set $K\subset \R^d$ a simplex,  ball or cube.
\end{abstract}

Suppose that $K\subset \R^d$ is the closure of its interior and compact. We note that the dimension of the polynomials of degree at most $n$ in $d$   variables is
\[{\rm dim}({\cal P}_n(\R^d))=N_n(=N):={n+d\choose d}.\]

For a basis $\{p_1,\cdots,p_N\}$ of  ${\cal P}_n(\R^d)$ and
$N$ points ${\bf x}_1,\cdots,{\bf x}_N$ in $K$ we may form the Vandermonde determinant
\[{\rm vdm}({\bf x}_1,\cdots,{\bf x}_N):={\rm det}([p_j({\bf x}_i)]_{1\le i,j\le N}).\]

In case the Vandermonde determinant is non-zero,  the problem of interpolation at these points by polynomials of degree at most $n$ is regular,  and we may,  in particular,  construct the fundamental Lagrange polynomials $\ell_i({\bf x})$ of degree $n$ with the property that
\[\ell_i({\bf x}_j)=\delta_{ij}.\]

The interpolation operator based on the points of $X,$
\[\pi_X\,:\, C(K) \to {\cal P}_n(\R^d); \quad \pi_X(f)=\sum_{i=1}^N f({\bf x}_i))\ell_i({\bf x})\]
has operator norm the Lebesgue constant defined as
\[ \Lambda_n(X;K):=\max_{{\bf x}\in K}\sum_{i=1}^N |\ell_i({\bf x})|.\]

In the case that $K\subset \R^d$ is a simplex, it was shown in \cite{B83} that, for $X$ the so-called equally spaced points (or simplex points) of degree $n$ there is a upper bound for the Lebesgue constant,  {\it independent} of the dimension. Specifically,
\[ \Lambda_n(X;K)\le {2n-1\choose n},\quad \forall d\ge1.\]
It follows that there is also such a bound for the Lebesgue points, i.e. those for which the Lebesgue constant is a minimum (and hence also highle likely for any good set of interpolation points).

The subject of this short paper is to study the growth {\it as a function of the dimension $d,$} of $\Lambda_n(X)$ for $n=1$ fixed and $X$ a set of degree one Fekete points for $K$ either a simplex,  ball or cube. We will see that for a ball or cube, the Lebesgue constant does grow with the dimension, making this property of the simplex somewhat unique and hence notable.

\begin{definition} A set $F\subset K$ of $N$ distinct points is said to a set of {\it Fekete} points of degree $n$ if they maximize 
$|{\rm vdm}({\bf x}_1,\cdots,{\bf x}_N)|$ over $K^N.$
\end{definition}

\begin{definition} A set  $F\subset K$ of $N$ distinct points is said to be a Fej\`er set if
\[\max_{x\in K} \sum_{i=1}^{N}\ell_i^2(x)=1.\]
\end{definition}

It is shown in \cite{B83} that a Fej\`er set is automatically also a Fekete set, but the reverse implication is not true. This problem is discussed in some detail in \cite{Bo}.

In the degree one case we have $N_1=d+1$ points  $X=\{{\bf x}_1,\cdots,{\bf x}_N\}\subset \R^d$ forming the vertices of a simplex. Writing
\[ {\bf x}_i=(x^{(i)}_1,x^{(i)}_2,\cdots, x^{(i)}_d)\in \R^d\]
the Vandermonde determinant for the basis of monomials $\{x_1,x_2,\cdots,x_d\}$ becomes
\[{\rm vdm}({\bf x}_1,\cdots,{\bf x}_{d})=\left| \begin{array}{cccccc}
1&x_1^{(1)}&x_2^{(1)}&\cdot&\cdot&x_d^{(1)}\cr
1&x_1^{(2)}&x_2^{(2)}&\cdot&\cdot&x_d^{(2)}\cr
\cdot&&&&&\cdot\cr
\cdot&&&&&\cdot\cr
1&x_1^{(d)}&x_2^{(d)}&\cdot&\cdot&x_d^{(d)}\cr
1&x_1^{(d+1)}&x_2^{(d+1)}&\cdot&\cdot&x_d^{(d+1)}\end{array}\right|  \]
which equals $\pm d!$ times the volume of the simplex with vertices the points of $X.$ Hence, in the degree one case, the Fekete points correspond to the vertices of the simplex of maximum volume contained in $K.$

 \medskip

\section{$K$ a Simplex}

In this case the Fekete points are just the vertices of the simplex and the associated Lagrange polynomials of degree one are linear polynomials which are one at a vertex and identically zero on the face opposite that vertex. Hence they are each non-negative on the simplex and
\[\sum_{i=1}^{d+1}|\ell_i(x)|=\sum_{i=1}^{d+1}\ell_i(x)=1,\quad x\in K.\]
In particular 
\[\Lambda_n(X;K)\equiv 1,\quad \forall d\ge 1.\]

\section{$K$ the Unit Ball}

We construct recusrsively the $d+1$ vertices of a regular simplex inscribed in the unit sphere in $\R^d$ as follows. Let
\[ X_1=\left[ \begin{array}{l} -1 \cr +1\end{array}\right]\in \R^{2\times 1}.\] 
Then, for $d>1$ and having defined $X_{d-1}\in R^{d\times (d-1)},$  let
\[ X_d:= \left[\begin{array}{cccc}
&&&-1/d\cr
&R_dX_{d-1}&&-1/d\cr
&&&-1/d\cr
0&\cdot&0&1\end{array}\right]\in \R^{(d+1)\times d} \]
where
\[R_d:=\frac{\sqrt{d^2-1}}{d}.\]
Each row of $X_d$ represents the cartesian coordinates of a point in $\R^d.$ Geometrically, $X_d$ is obtained by placing the points of $X_{d-1},$ properly scaled, on the ball in $\R^{d-1}$ given by the intersection of the ball in $\R^d$ with the level set $x_d=-1/d$  and the adding the "north pole" $(0,\cdots,0,1)\in \R^d.$ We claim that these form the verticess of a regular simplex inscribed in the unit sphere. This will be evident from the following lemmas.

\begin{lemma} \label{lem1}The centroid of the points gievn by $X_d$ is ${\bm 0}_d\in \R^d,$ i.e., $X_d^t{\mathbbm{1}}_{d+1}={\bm 0}_d.$
\end{lemma}
\noindent {\bf Proof} (by induction on $d$). For $d=1,$
\[X_1^t\mathbbm{1}_2=[-1\,\,1]\left[\begin{array}{c}1\cr 1\end{array}\right]=0.\]
Hence assume that $X_{d-1}^t\mathbbm{1}_d ={\bm 0}_{d-1}$ and calculate
\begin{align*}
X_d^t\mathbbm{1}_{d+1}&= \left[\begin{array}{cccc}
&&&0\cr
&R_dX_{d-1}^t&&0\cr
&&&0\cr
-1/d&\cdot&-1/d&1\end{array}\right]\,\left[\begin{array}{c}1\cr \cdot\cr\cdot\cr 1\end{array}\right]\cr
&=\left[ \begin{array}{c}R_dX_{d-1}^t\mathbbm{1}_d\cr d(-1/d)+1\end{array}\right]\cr
&={\bm 0}_d.
\end{align*}
$\square$
\bigskip
\begin{lemma} \label{lem2} We have
\[ X_d^tX_d=\frac{d+1}{d} I_d.\]
\end{lemma}
\noindent {\bf Proof} (by induction on the dimension $d$). For $d=1,$
\[ X_1^tX_1=[-1 \,\,1]\left[\begin{array}{c}-1\cr+1\end{array}\right]=2=\frac{1+1}{1}I_1.\]
Assume then the Lemma holds for dimension $d-1.$ Then
\begin{align*}
X_d^tX_d&=\left[\begin{array}{cccc}
&&&0\cr
&R_dX_{d-1}^t&&0\cr
&&&0\cr
-1/d&\cdot&-1/d&1\end{array}\right]
\,\left[\begin{array}{cccc}
&&&-1/d\cr
&R_dX_{d-1}&&-1/d\cr
&&&-1/d\cr
0&\cdot&0&1\end{array}\right]\cr
&=\left[\begin{array}{cc}
R_d^2X_{d-1}^tX_{d-1}& R_dX_{d-1}^t\left[\begin{array}{c}-1/d\cr\cdot\cr\cdot\cr-1/d\end{array}\right]\cr
\bigl[-1/d\,\cdot\,\cdot\,-1/d\bigr] R_dX_{d-1}&\frac{d}{d^2}+1\end{array}\right]\cr
&=\left[ \begin{array}{cc}\frac{d^2-1}{d^2}\,\frac{d}{d-1}I_{d-1}& {\bm 0}_{d-1}\cr
{\bm 0}^t_{d-1}& 1+1/d\end{array}\right]\cr
&=\frac{d+1}{d}I_d
\end{align*}
as
\[ X_{d-1}^t\left[\begin{array}{c}-1/d\cr\cdot\cr\cdot\cr-1/d\end{array}\right]=-\frac{1}{d}
X_{d-1}^t\left[\begin{array}{c}1\cr\cdot\cr\cdot\cr 1\end{array}\right]={\bm 0}_{d-1}\]
by Lemma \ref{lem1}. $\square$

\bigskip
\begin{lemma} \label{lem3} We have
\[ X_dX_d^t=\frac{d+1}{d} I_{d+1}-\frac{1}{d} \mathbbm{1}_{d+1} \mathbbm{1}_{d+1}^t\in \R^{(d+1)\times (d+1)}.\]
\end{lemma}
\noindent{\bf Proof} (by induction on the dimension $d$). 
For $d=1,$
\begin{align*}
X_1 X_1^t &= \left[\begin{array}{c}-1\cr+1\end{array}\right]\, [-1\,\,+1]=\left[\begin{array}{cc}
+1&-1\cr -1&+1\end{array}\right]\cr
&=\left[\begin{array}{cc}2&0\cr 0&2\end{array}\right]-\left[\begin{array}{cc}
1&1\cr1&1\end{array}\right]\cr
&=\frac{1+1}{1}I_2-\frac{1}{1} \mathbbm{1}_{2} \mathbbm{1}_{2}^t.
\end{align*}
Hence assume that the Lemma holds for $d-1.$ Then
\begin{align*}
X_d X_d^t&=
\left[\begin{array}{cccc}
&&&-1/d\cr
&R_dX_{d-1}&&-1/d\cr
&&&-1/d\cr
0&\cdot&0&1\end{array}\right]\,
\left[\begin{array}{cccc}
&&&0\cr
&R_dX_{d-1}^t&&0\cr
&&&0\cr
-1/d&\cdot&-1/d&1\end{array}\right]\cr
&=\left[\begin{array}{cc}
R_d^2X_{d-1}X_{d-1}^t+\frac{1}{d^2} \mathbbm{1}_d\mathbbm{1}_d^t & 
\begin{array}{c}-1/d\cr\cdot\cr\cdot\cr -1/d\end{array}\cr
\begin{array}{cccc}-1/d&\cdot&\cdot&-1/d\end{array}& 1\end{array}\right]\cr
&=\left[\begin{array}{cc}
\frac{d^2-1}{d^2}\Bigl(\frac{d}{d-1}I_d -\frac{1}{d-1} \mathbbm{1}_d\mathbbm{1}_d^t\Bigr)+\frac{1}{d^2} \mathbbm{1}_d\mathbbm{1}_d^t & 
\begin{array}{c}-1/d\cr\cdot\cr\cdot\cr -1/d\end{array}\cr
\begin{array}{cccc}-1/d&\cdot&\cdot&-1/d\end{array}& 1\end{array}\right]\cr
&=\left[\begin{array}{cc}
\frac{d+1}{d}I_d -\frac{1}{d} \mathbbm{1}_d\mathbbm{1}_d^t & 
\begin{array}{c}-1/d\cr\cdot\cr\cdot\cr -1/d\end{array}\cr
\begin{array}{cccc}-1/d&\cdot&\cdot&-1/d\end{array}& 1\end{array}\right]\cr
&=\frac{d+1}{d} I_{d+1} -\frac{1}{d} \mathbbm{1}_{d+1}\mathbbm{1}_{d+1}^t.
\end{align*}
$\square$

\bigskip \noindent {\bf Remark}. The matrix $X_dX_d^t$ is the Gram matrix of the inner products of the points (rows) of $X_d.$ The diagonal entries are all $1$ indicating that each point is on the unit sphere. The off-diagonal entries are all $-1/d$ indicating that the points are equi-angular. The square of the distance between ${\bf x}_i$ and ${\bf x}_j$ is
\begin{align*}
\|{\bf x}_i-{\bf x}_j\|_2^2&= ({\bf x}_i - {\bf x}_j)^t({\bf x}_i-{\bf x}_j)\cr
&={\bf x}_i^t{\bf x}_j-2{\bf x}_i{\bf x}_j+{\bf x}_j^t{\bf x}_j\cr
&=1-2(-1/d)+1=2\frac{d+1}{d},\quad \forall i\neq j.
\end{align*}
$\square$

\bigskip The barycentric coordinates of ${\bf y}\in \R^d$ with respect to the points of $X_d$ are given by the vector ${\bm \lambda}\in \R^{d+1}$ such that\par
\noindent (a) $X_d^t{\bm \lambda}={\bf y},$ and\par
\noindent (b) $\mathbbm{1}_{d+1}^t{\bm \lambda}=1.$
\begin{lemma}\label{lem4}
We have, for ${\bf y}\in \R^d,$
\[ {\bm \lambda}=\frac{1}{d+1}\mathbbm{1}_{d+1}+\frac{d}{d+1} X_d{\bf y}.\]
\end{lemma}
\noindent {\bf Proof}. We need just to verify the properties (a) and (b). For (a),
\begin{align*}
X_d^t{\bm \lambda}&=X_d^t\left\{\frac{1}{d+1}\mathbbm{1}_{d+1}+\frac{d}{d+1} X_d{\bf y}\right\}\cr
&=\frac{1}{d+1} X_d^t\mathbbm{1}_{d+1}+\frac{d}{d+1}(X_d^tX_d){\bf y}\cr
&= {\bm 0}_{d+1}+\frac{d}{d+1}\left(\frac{d+1}{d}I_d\right){\bf y}\quad(\hbox{by\,\,Lemmas\,\,1\,\,and\,\,2})\cr
&={\bf y}.
\end{align*}
For (b), 
\begin{align*}
\mathbbm{1}_{d+1}^t{\bm \lambda}&=\mathbbm{1}_{d+1}^t\left\{  \frac{1}{d+1}\mathbbm{1}_{d+1}+\frac{d}{d+1} X_d{\bf y} \right\}\cr
&=\frac{1}{d+1} \bigl( \mathbbm{1}_{d+1}^t \mathbbm{1}_{d+1}\bigr)+\frac{d}{d+1}\bigl( X_d^t \mathbbm{1}_{d+1}\bigr)^t{\bf y}\cr
&=\frac{1}{d+1}(d+1)+0\quad (\hbox{by\,\,Lemma\,\,1})\cr
&=1.
\end{align*}
$\square$

\bigskip

It is easy to confirm that the Lagrange polynomials (of degree one) are precisely the barycentric coordinates. Hence the Lebesgue function is
\[ \Lambda_1({\bf y})=\sum_{j=1}^{d+1} |\lambda_j({\bf y})|=\|{\bm \lambda}({\bf y})\|_1\]
and the Lebesgue constant is
\[ \Lambda_1=\max_{\|{\bf y}\|_2\le 1}\Lambda_1({\bf y})=\max_{\|{\bf y}\|_2\le 1} \|{\bm \lambda}({\bf y})\|_1.\]
\begin{lemma}\label{lem4} For 
\[ {\bm \lambda}=\frac{1}{d+1}\mathbbm{1}_{d+1}+\frac{d}{d+1} X_d{\bf y}\]
we have
\[\sum_{j=1}^{d+1}\lambda_j^2 = \frac{1+d\|{\bf y}\|_2^2}{d+1}\le 1,\,\,\|{\bf y}\|_2\le1.\]
\end{lemma}
\noindent {\bf Proof}. We calculate
\begin{align*}
\sum_{j=1}^{d+1} \lambda_j^2&= {\bm \lambda}^t{\bm \lambda}\cr
&= \frac{1}{(d+1)^2}\bigl\{ \mathbbm{1}_{d+1}^t+d{\bf y }^tX_d^t\bigr\}\,
\bigl\{\mathbbm{1}_{d+1}+d X_d {\bf y}\bigr\}\cr
&=  \frac{1}{(d+1)^2}\bigl\{ \mathbbm{1}_{d+1}^t \mathbbm{1}_{d+1}+2d ( \mathbbm{1}_{d+1}^tX_d){\bf y} +d^2{\bf y}^t(X_d^tX_d){\bf y}\bigr\}\cr
&=\frac{1}{(d+1)^2}\bigl\{  (d+1)+0+d^2{\bf y}^t\bigl( \frac{d+1}{d} I_d\bigr){\bf y} \bigl\}\cr
&= \frac{1}{(d+1)^2}\bigl\{ (d+1)+d(d+1)\|{\bf y}\|_2^2\bigr\}\cr
&=\frac{1+d\|{\bf y}\|_2^2}{d+1}.
\end{align*}
$\square$

\bigskip
It follows that the points of $X_d$ form a F\'ejer and hence Fekete set and, in particular, 
\[\Lambda_1\le \sqrt{d+1}.\]
We claim that, in fact, this upper bound is the correct order of growth.

\begin{prop} \label{prop1} For $K$ the unit ball, we have
\[ \sqrt{d} \le \Lambda_1 \le \sqrt{d+1}.\]
\end{prop}
\noindent {\bf Proof}. We may write
\begin{align*}
\Lambda_1&= \max_{\|{\bf y}\|_2\le 1} \|{\bm \lambda}({\bf y})\|_1\cr
&=  \max_{\|{\bf y}\|_2\le 1}\,\, \max_{{\bm \epsilon}\in\{\pm 1\}^{d+1}} {\bm \epsilon}^t{\bm \lambda}({\bf y})\cr
&= \max_{{\bm \epsilon}\in\{\pm 1\}^{d+1}}\,\, \max_{\|{\bf y}\|_2\le 1} {\bm \epsilon}^t{\bm \lambda}({\bf y})\cr
&= \max_{{\bm \epsilon}\in\{\pm 1\}^{d+1}}\,\, \max_{\|{\bf y}\|_2\le 1} {\bm \epsilon}^t
\bigl\{ \frac{1}{d+1}\mathbbm{1}_{d+1}+\frac{d}{d+1} X_d {\bf y}\bigr\}\cr
&=\max_{{\bm \epsilon}\in\{\pm 1\}^{d+1}}\,\, \max_{\|{\bf y}\|_2\le 1} 
\bigl\{ \frac{s}{d+1}+\frac{d}{d+1}({\bm \epsilon}^t X_d) {\bf y}\bigr\}\cr
&\qquad \quad({\rm with}\,\, s:= \sum_{j=1}^{d+1} \epsilon_j)\cr
&=\max_{{\bm \epsilon}\in\{\pm 1\}^{d+1}} \frac{1}{d+1}\bigl\{
s+d\| {\bm \epsilon}^tX_d\|_2\bigr\}\cr
&\qquad\quad ({\rm with}\,\, {\bf y}=({\bm \epsilon}^tX_d)^t/\| {\bm \epsilon}^tX_d\|_2)\cr
&= \max_{{\bm \epsilon}\in\{\pm 1\}^{d+1}} \frac{1}{d+1}\bigl\{
s+d\sqrt{{\bm \epsilon}^t(X_dX_d^t){\bm \epsilon}}\bigr\}\cr
&=\max_{{\bm \epsilon}\in\{\pm 1\}^{d+1}} \frac{1}{d+1}\left\{
s+d\sqrt{{\bm \epsilon}^t\Bigl\{ \frac{d+1}{d} I_{d+1}-\frac{1}{d} \mathbbm{1}_{d+1}\mathbbm{1}_{d+1}^t\Bigr\}{\bm \epsilon}}\right\}\cr
&\qquad\quad ({\rm by\,\, Lemma\,\, 3})\cr
&= \max_{{\bm \epsilon}\in\{\pm 1\}^{d+1}} \frac{1}{d+1}\left\{
s+d\sqrt{ \frac{d+1}{d} \|{\bm \epsilon}\|_2^2-\frac{s^2}{d}}\right\}\cr
&=\max_{{\bm \epsilon}\in\{\pm 1\}^{d+1}} \frac{1}{d+1}\left\{
s+\sqrt{d}\sqrt{(d+1)^2-s^2}\right\}\cr
&\qquad\quad ({\rm as}\,\,\|{\bm \epsilon}\|_2^2=d+1)\cr
\end{align*}
where again
\[s:=\sum_{j=1}^{d+1}\epsilon_j \in\{-(d+1),\cdots,0,\cdots,(d+1)\}.\]
However, if  $s\le0,$ we may replace ${\bm \epsilon}$ by $-{\bm \epsilon}$ and hence we may assume that $s\ge0,$ i.e., $s\in\{0,1,\cdots,(d+1)\}.$

For the continuous function
\begin{align*}
f(s)&:=\frac{1}{d+1}\bigl\{s+\sqrt{d}\sqrt{(d+1)^2-s^2}\bigr\}, \,\,s\in[0,d+1],\cr
f'(s)&=\frac{1}{d+1}\bigl\{1-\sqrt{d}\frac{s}{\sqrt{(d+1)^2-s^2}}\bigr\}
\end{align*}
and has a single maximum at the critical point given by 
\begin{align*}
\sqrt{(d+1)^2-s^2}&=\sqrt{d}\,s\cr
\iff (d+1)^2-s^2&=d s^2\cr
\iff s&=\sqrt{d+1}
\end{align*}
for which 
\[f(s)=\sqrt{d+1}.\]
Hence, as already noted
\[\Lambda_1\le \sqrt{d+1}.\]
However, $s=\sqrt{d+1}$ is rarely an integer and hence this upper bound is only attained in special circumstances.

Now note that for $d$ {\it odd}, $s=0$ is attainable (half the $\epsilon_j=+1$ and the other half equal to $-1$) and hence
\[\Lambda_1\ge \frac{1}{d+1}\bigl\{0+\sqrt{d}\sqrt{(d+1)^2-0}\bigr\}=\sqrt{d}.\]
It follows that, for $d$ odd,
\[\sqrt{d}\le \Lambda_1\le \sqrt{d+1}.\]
In case $d$ is even then $s=1$ is attainable and hence, for $d$ {\it even}
\begin{align*}
\Lambda_1&\ge \frac{1}{d+1}\bigl\{1+\sqrt{d}\sqrt{(d+1)^2-1}\bigr\}\cr
&=\frac{1+d\sqrt{d+2}}{d+1}\cr
&\ge \sqrt{d}
\end{align*}
as is easily confirmed. $\square$

\section{The Case of $K=[-1,1]^d$ a Cube}

Since the Vandermonde determinant is linear as a function of each point separately,  its maximum will be attained at a subset of the vertices of the cube $[-1,1]^d,$ i.e., the optimal Vandermonde matrix is a matrix of all $\pm 1$ entries.  To determine the maximum determinant of such $\pm1$ matrices is the celebrated Hadamard Determinant Problem (1893) ,  whose solution is not yet known in general dimension.   There are however special dimensions in which the solution of the Hadamard problem is known and given by a so-called Hadamard matrix where the rows and columns are mutually orthogonal.  
Correspondingly,  in dimensions $d$ for which there exists a Hadamard matrix of order $n=d+1$ the
$d+1$ Fekete points can be explicitly expressed in terms of the rows of the Hadamard matrix.

\begin{definition} A matrix $H\in\R^{n\times n}$ with entries $H_{ij}\in{\pm1}$ and rows and columns  orthogonal,  i.e.,
\[H_nH_n^t=nI_n\]
is said to be a Hadamard matrix.
\end{definition}

Sylvester's construction gives a Hadamard matrix for all $n$ a power of 2,  but the existence for many other values of $n$ is also known.

Now suppose that $d$ is such that a Hadmard matrix,  $H_{d+1},$ of dimension $d+1$ exists. By multiplying on the left and right by  appropriate diagonal matrices, we may assume that the first row and first column of $H_{d+1}$ are all $1s.$ We let $X_d\in \R^{(d+1)\times d}$ be the matrix obtained by removing the first column of $H_{d+1}.$ The $d+1$ rows of $X_d$ give the coordinates of a subset of $d+1$ vertices of the cube $[-1,1]^d,$ and it is these points that we consider.  In particular $V_d:=H_{d+1}$ is the Vandermonde matrix for these points and the polynomials of degree at most one with basis
\[\{1,x_1,\cdots,x_d\}.\]
Hence,  by the definition of Hadamard matrices the points $X$ are such that their associated Vandermonde matrix has determinant
as large as possible (in absolute value) and hence are Fekete points. 

\medskip 
\noindent {\bf Example}.  For $d=3,$ 
\[H_{d+1}=H_4=\left[\begin{array}{rrrr}
1 &1 & 1  &1 \cr
1 &-1 &1 &-1 \cr
1 &1 & -1& -1 \cr
1&-1&-1&1 \cr
\end{array} \right] \]
 so that the four points are
 \[(1  ,   1  ,   1), \,
    (-1   ,  1  ,  -1),\,
     (1   , -1  ,  -1),\,
    (-1  ,  -1  ,   1).\]
The simplex with these vertices is shown in Figure 1 below.

\begin{center}\begin{figure}[t!]\centering
 \includegraphics[width=.6\textwidth,height=.6\textheight]{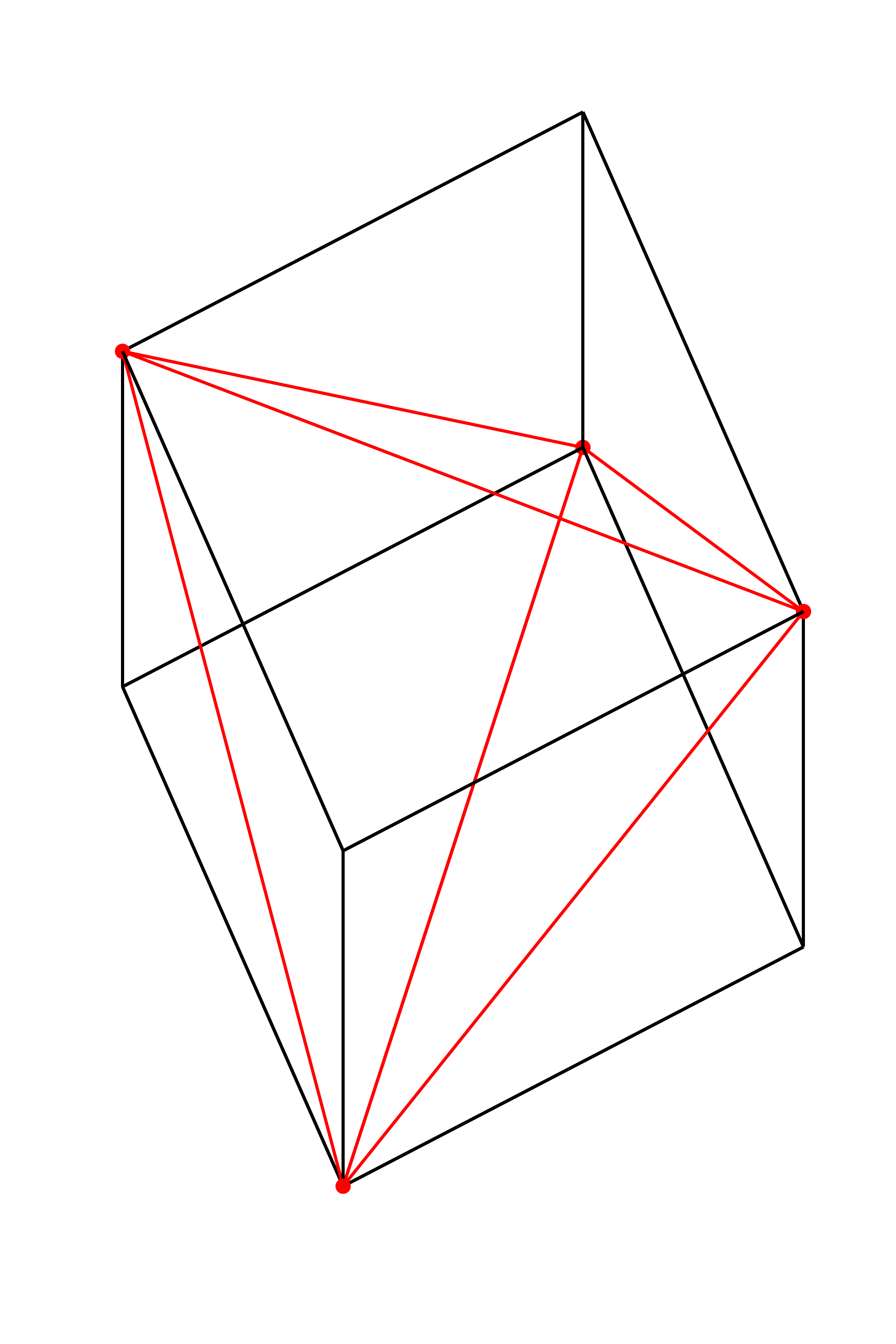}\hfill
 \caption{Regular Simplex Inscribed in the Cube}
\end{figure}
\end{center}

The associated fundamental Lagrange polynomials are
\[[\ell_1({\bf x}),\cdots,\ell_{d+1}({\bf x})]=[1, {\bf x}^t]V_d^{-1}
=\frac{1}{d+1}[1, {\bf x}^t]V_d^{t}.\]
They have the property that
\begin{eqnarray*}
\sum_{i=1}^{d+1}\ell_i^2({\bf x})&=&[\ell_1({\bf x}),\cdots,\ell_{d+1}({\bf x})]\times \left[\begin{array}{c}
\ell_1({\bf x})\cr
\cdot\cr
\cdot\cr
\ell_{d+1}({\bf x})
\end{array}\right]\cr
&=&\frac{1}{(d+1)^2}[1, {\bf x}^t]V_d^{t}V_d\left[\begin{array}{c} 1\cr {\bf x}\end{array}\right]\cr
&=&\frac{1+\|{\bf x}\|_2^2}{d+1}\le 1
\end{eqnarray*}
for all ${\b x}\in[-1,1]^d.$

In other words,  they are also a set of F\'ejer points.  We mention again that,  as shown in \cite{B83}, this is also sufficient to prove that the points $X_d$ are Fekete points.  

\noindent {\bf Remark}.  Such points $X_d$ form the vertices of a {\it  regular} simplex.  As the Vandermonde determinant is a (dimensional) multiple of the volume of this simplex,  it is of maximal volume.  Also,  as the sum of the Lagrange polynomials squared is bounded by 1 on the circumball $B_d:={{\bf x}\in \R^d\,:\, \|{\bf x}\|_2\le \sqrt{d}},$ $X_d$ is also a Fekete set  for $B_d.$ $\square$

In particular,  we have again that
\[\Lambda_1\le\sqrt{d+1}.\]

We claim that in certain dimensions this upper bound is also attained.

\begin{prop} In case $d=m^2-1$ is such that the Hadamard matrices $H_m$ and $H_{d+1}=H_m \otimes  H_m$ exist  then 
\[\Lambda_1=\sqrt{d+1}.\]
\end{prop}

\noindent {\bf Proof}.  We have
\begin{align*}
\Lambda_1 &= \max_{{\bf x}\in [-1,1]^d} \sum_{j=1}^{d+1}|\ell_j({\bf x})|\cr
&=\max_{{\bf x}\in \{\pm 1\}^d}\sum_{j=1}^{d+1}|\ell_j({\bf x})|\cr
&= \max_{{\bf x}\in \{\pm 1\}^d}\max_{{\bm \epsilon}\in\{\pm1\}^{d+1}}\sum_{j=1}^{d+1}\epsilon_j\ell_j({\bf x})\cr
&=\max_{{\bf x}\in \{\pm 1\}^d}\max_{{\bm \epsilon}\in\{\pm1\}^{d+1}} \frac{1}{d+1}[1\,\,{\bf x}^t]H_{d+1}^t{\bm \epsilon}.
\end{align*}
Hence it suffices to exhibit ${\bf x}\in\{\pm1\}^d,$ ${\bm \epsilon}\in\{\pm 1\}^{d+1}$ such that
\[\frac{1}{d+1}[1\,\,{\bf x}^t]H_{d+1}^t{\bm \epsilon}=\sqrt{d+1}.\]
Now write,  in columns,
\[ H_m =[{\bf h}_1\,\,{\bf h}_2,\cdots,{\bf h}_m],\]
with ${\bf h}_j\in\R^m,$
and set
\[{\bm \epsilon}=\left[\begin{array}{c}
{\bf h}_1\cr {\bf h}_2\cr \cdot\cr\cdot\cr {\bf h}_m\end{array}\right]\in\R^{m^2}.\]

Note that $H_m^tH_m=mI_m$ means that
\[ H_m^t {\bf h}_j =m {\bf e}_j\]
the canonical basis vector.

Further,  as $H_{d+1}=H_{m^2}=H_m \otimes H_m,$ we may write, in block form,
\[H_{d+1}^t=[H_{ij}]_{1\le i,j\le m},\,\,\,H_{ij}=\pm H_m^t.\]
It follows that the $i$th block of $H_{d+1}^t{\bm \epsilon}$ is
\begin{align*}
(H_{d+1}^t{\bm \epsilon})_i&= \sum_{j=1}^m (\pm 1)H_m{^t\bf h}_j\cr
&= \sum_{j=1}^m (\pm m){\bf e}_j
\end{align*}
which is a vector all of whose components are $\pm m.$ In particular $H_{d+1}^t{\bm \epsilon}\in \R^{m^2}$ is also a vector with all components $\pm m.$ Without loss of generality we may assume that the first component is $+1.$ Taking then
\[x_j:= {\rm sgn}(H_{d+1}{\bm \epsilon})_{j+1}\]
we obtain
\begin{align*}
\frac{1}{d+1}[1\,\,{\bf x}^t]H_{d+1}^t{\bm \epsilon}&=\frac{1}{d+1}
\| H_{d+1}^t{\bm \epsilon}\|_1\cr
&=\frac{1}{d+1}(d+1)m\cr
&=m=\sqrt{d+1}.
\end{align*}
$\square$

\section{$K$ the Complex Torus}

There are analogous results for the complex version of the cube, the Torus. Consider
\[K=\mathbb{T}^d:=\{{\bf z}\in \C^d\,:\, |z_j|=1,\,\,1\le j\le d\}.\]
In this case the classical Fourier matrix plays the role of the Hadamard matrix.

\begin{definition} The Fourier matrix $F_n\in \C^{n\times n}$ defined by
\[ F_n:=[\omega^{jk}]_{1\le j,k\le n},\quad \omega:=\exp(2\pi i/n)\]
is known as the Fourier matrix.
\end{definition}

As is well known,  the Fourier matrix has orthogonal rows and columns, i.e.,
\[F_n^*F_n=nI_n\]
and is sometimes referred to as a complex Hadamard matrix,  as the entries all have modulus $1.$

Just as for the cube and Hadamard matrix we let
$X_d\in \C^{(d+1)\times d}$ be the matrix obtained by removing the first column of $F_{d+1}.$ The $d+1$ rows of $X_d$ give the coordinates of a subset of $d+1$ points of the torus $\mathbb{T}^d.$ In particular $V_d:=F_{d+1}$ is the Vandermonde matrix for these points and the polynomials of degree at most one with basis
\[\{1,z_1,\cdots,z_d\}.\]

The associated fundamental Lagrange polynomials are
\begin{align*}
\left[\begin{array}{c}
\ell_1({\bf z})\cr\cdot\cr\cdot\cr \ell_{d+1}({\bf z})\end{array}
\right]&=V_d^{-t}\left[\begin{array}{c}1\cr {\bf z}\end{array}\right]\cr
&=V_d^{-1}\left[\begin{array}{c}1\cr {\bf z}\end{array}\right]
\quad ({\rm as}\,\,F_n^t=F_n)\cr
&=\frac{1}{d+1}V_d^{*}\left[\begin{array}{c}1\cr {\bf z}\end{array}\right].
\end{align*}
They have the property that
\begin{eqnarray*}
\sum_{i=1}^{d+1}|\ell_i({\bf z})|^2&=&\left[\begin{array}{c}
\ell_1({\bf z})\cr
\cdot\cr
\cdot\cr
\ell_{d+1}({\bf z})
\end{array}\right]^*\times \left[\begin{array}{c}
\ell_1({\bf z})\cr
\cdot\cr
\cdot\cr
\ell_{d+1}({\bf z})
\end{array}\right]\cr
&=&\frac{1}{(d+1)^2}[1, {\bf z}^*]V_dV_d^*\left[\begin{array}{c} 1\cr {\bf z}\end{array}\right]\cr
&=&\frac{1+\|{\bf z}\|_2^2}{d+1}\le 1
\end{eqnarray*}
for all ${\bf z}\in \mathbb{T}^d.$

In other words,  they are also a set of F\'ejer points. and this is  sufficient to show that the points $X_d$ are Fekete points.  

\begin{prop} For $K=\mathbb{T}^d, $ the complex torus,  and $d=m^2-1$ for any positive integer $m,$
\[\Lambda_1=\sqrt{d+1}.\]
\end{prop}
\noindent {\bf Proof}.  From the fact that $X_d$ are F\'ejer points we have the upper bound,
\[\Lambda_1\le \sqrt{d+1}.\]
To show the lower bound we argue as for the real cube with $H_{d+1}$ replaced by $F_{d+1},$ using the fact that for $d+1=m^2,$ $F_{d+1}=F_m\otimes F_m.$ In particular
\begin{align*}
\Lambda_1 &= \max_{{\bf z}\in \mathbb{T}^d} \sum_{j=1}^{d+1}|\ell_j({\bf z})|\cr
&= \max_{{\bf z}\in \mathbb{T}^d}\max_{|\epsilon_j|\le1,\,1\le j\le(d+1)}\sum_{j=1}^{d+1}\epsilon_j\overline{\ell_j({\bf z})}\cr
&=\max_{{\bf z}\in \mathbb{T}^d}\max_{|\epsilon_j|\le1,\,1\le j\le(d+1)} \frac{1}{d+1}[1\,\,{\bf z}^*]F_{d+1}{\bm \epsilon}.
\end{align*}
Hence it suffices to exhibit ${\bf z}\in\mathbb{T}^d$ and ${\bm \epsilon}\in \C^{d+1}$  with $|\epsilon_j|\le1, $ $1\le j\le d+1,$ such that
\[\frac{1}{d+1}[1\,\,{\bf z}^t]F_{d+1}{\bm \epsilon}=\sqrt{d+1}.\]
It is easy to verify that ${\bm \epsilon}$ the stacked columns of $F_m^*$ and ${\bf z}$ a suitably cjosen complex sign vector, have this property. $\square$

\section*{Acknowledgements}
 RITA 
``Research ITalian network on Approximation''.

\end{document}